\documentclass{article}
\usepackage{amsmath,amsfonts,amssymb}

\numberwithin{equation}{section}

\newcommand{\sezione}[1]{\section{\normalsize #1}}

\newtheorem{theorem}[equation]{\textnormal{THEOREM}}
\newtheorem{lemma}[equation]{\textnormal{LEMMA}}
\newtheorem{corollary}[equation]{\textnormal{COROLLARY}}

\newtheorem{definition}[equation]{\textnormal{DEFINITION}}
\newtheorem{example}[equation]{\textnormal{EXAMPLE}}

\def\NN{\Bbb N}
\def\QQ{\Bbb Q}
\def\ZZ{\Bbb Z}

\def\spec{\mbox{Spec }}

\begin{document}

\vspace*{2.8cm}
\section*{An introduction to tight closure}
\begin{flushleft}
\medskip
KAREN E. SMITH Department of Mathemathics,
University of Michigan,
Ann Arbor, MI, 48109, kesmith@math.lsa.umich.edu\par
\end{flushleft}
\bigskip
\bigskip

\vspace{2\baselineskip}
\sezione{INTRODUCTION}

This is an expanded version of my lecture at the 
Conference in Commutative Algebra and Algebraic Geometry in  Messina Italy in June 1999.
The purpose of the talk was to  give a brief introduction to the subject of
tight closure, aimed at  commutative algebraists who have not before studied this topic.
 The first part focused mainly on the definition and basic properties,  with the
second part focusing  
on some applications to algebraic geometry, particularly to global generation of adjoint linear
series.   These lecture notes follow even more closely a series of two lectures I gave in
Kashikojima, Japan, at the Twentieth Annual Japanese Conference in Commutative
Algebra the previous fall, and were distributed also in conjunction with that conference.
 I wish to thank the organizers of both conferences,
Professors Restuccia and Herzog for the European conference, and 
 Professors Hashimoto and Yoshida, for the Japanese conference. Both events
were a smashing success. Special thanks are due also to Rosanna Utano,
for help editing the tex file.

\bigskip
Tight closure
was introduced by Mel Hochster and Craig Huneke
 in 1986 \cite{HH1}.
Today it is still a subject of very active research, with an ever increasing list 
of applications.  
Applications include areas like the study of Cohen-Macaulayness. For example, the famous
Hochster-Roberts theorem on the Cohen-Macaulayness of rings of invariants has a
simple tight closure proof \cite{HH1}.  Also, the existence of big Cohen-Macaulay algebras for
rings containing a field  was proved with ideas from tight closure \cite{HH3}, and  the
existence of ''arithmetic  Macaulayfications"  in some cases was discovered with tight closure
\cite{AHS},
\cite{Ku}. Tight closure has provided us with  greater insight into integral closure, and into 
the homological theorems that
grew out of Serre's work on multiplicities.  For example,  it gives us simple proofs of the
Brian\c con-Skoda Theorem, the 
Syzygy Theorem of Evans and Griffith  and  of the monomial conjecture (in
equi-characteristic) \cite{HH1}.  Tight closure provided the 
inspiration for  results on the simplicity of rings of differential operators on certain
rings of invariants \cite{SV}, and it has produced "uniform" Artin-Rees theorems \cite{Hu1}.
There are also numerous applications  to and connections with algebraic geometry,  such as in
the study of singularities
\cite{W}, \cite{S1}, \cite{Ha1}, of vanishing theorems \cite{Ha2}, \cite{HS}, \cite{Ha4}, and of
adjoint linear series \cite{S4}, \cite{S7}.
  In Section 3, I will summarize some of these
applications to algebraic geometry, although of course, there will not be enough time to do any
of them any justice.

\bigskip
Let us begin with our first  task: to introduce the definition of tight closure
before tackling its main properties in the next section. 

\smallskip
 Tight closure is a closure operation performed on ideals in a commutative,
Noetherian ring containing a field (that is, of "equi-characteristic"). 
The tight closure of an ideal $I$ is an ideal $I^*$ containing $I$.
The definition is based on reduction to characteristic $p$, where the
Frobenius (or $p-th$ power map) is then used.  To keep things as
simple as possible, we treat only the characteristic $p$ case here.
\par
\bigskip

\begin{definition}
{\rm Let $R$  be a Noetherian domain of prime characteristic $p$, and let $I$ be an ideal
with generators $(y_1, \dots, y_r)$.  An element $z$ is defined to be in the tight closure
$I^*$ if  there exists a non--zero element $c$ of $R$ such that 
$$(*){\phantom{...............................}}cz^{p^e} \in (y_1^{p^e},
\dots, y_r^{p^e})
{\phantom{...............................}}$$ for all $e \gg 0$.}
\end{definition}
\par
\bigskip

Loosely speaking, the tight closure consists of all elements that are
"almost" in $I$  as far as the Frobenius map is concerned. Indeed, 
if we take the $p^e-th$ root of $(*)$ above,
we see that
$$
c^{1/p^e}z \in IR^{1/p^e}
$$
for $e \gg 0$.
As $e$ goes to infinity, $1/p^e$ goes to zero, so in some sense $c^{1/p^e}$
goes to 1 (this idea  can be made precise with valuations).
 So $z$ is "almost in"  $I$, at least after applying the Frobenius map.

It is not important  to restrict to the case where $R$ is a domain; we can define
tight closure in an arbitrary Noetherian ring of characteristic $p$ by requiring that $c$ is
not in  any minimal prime. However,  because most theorems about tight closure reduce to
 the domain 
case, we  treat only the domain case in this lecture.
\par
\bigskip

\begin{example}
{\bf Example.} Let $R$  be the hypersurface ring 
$$\frac{k[x, y, z]}{(x^3 + y^3 - z^3)},$$
where $k$ is any field whose characteristic is not 3.
Then 
$$
(x, y)^* = (x, y, z^2).
$$
Indeed,  if $k$ has characteristic $p$, we can write 
$$
(z^2)^q = (x^3 + y^3)^{\frac{2q-r}{3}}z^r  
$$
where $r = 1 $ or $2$ and $q = p^e$.  
Expanding this expression 
as
$$
z^r \sum \left( \begin{array}{c}
\frac{2q-r}{3}\\
i
\end{array}  \right) x^{3i}y^{2q-r -3i},
$$
it is easy to see that each monomial $x^m y^n$ appearing in the sum has 
either $m \geq q$ or $n \geq q$ unless both $m$ and $n$ equal $q -1$ (which only
happens in the case where $q = 1\ \ \mbox{ mod  $3$}$).
So we can take $c = x $ (or $y$), and conclude that
$c(z^2)^q \in (x^q, y^q)$ for all
$q = p^e$.
Thus
$z^2
\in (x, y)^*$. A similar argument can be used to show that $z$ is not in $(x, y)^*$.
Because this works for all  $p$ (except $p = 3$), we declare that 
$z^2 $, but not $z$,
 is in the tight closure of $(x, y)$ also in characteristic zero. 
So $(x, y)^* = (x, y, z^2)$ in every characteristic $p \geq 0$ except $p =3$.
\end{example}
\bigskip
\sezione{MAIN PROPERTIES OF TIGHT CLOSURE}
\bigskip
The definition of tight closure takes some getting used to.
Fortunately,  one can   understand many applications of  tight 
closure if one simply accepts the following properties of
tight closure as axioms:
\par
\bigskip
\noindent
{\bf Main Properties of Tight Closure}
{\it
\begin{enumerate}
\item[\em (1)] If $R$ is regular, then all ideals of $R$ are tightly closed.
\item[\em (2)] If $R \hookrightarrow S$ is an integral extension, then $IS\cap R \subset I^*$
for all ideals $I$ of $R$.
\item[\em (3)] If $R$ is local, with system of parameters $x_1, \dots, x_d$, then 
$(x_1, \dots, x_i):x_{i+1} \subset (x_1, \dots, x_i)^* $ 
("Colon Capturing").
\item[\em (4)]  If $\mu$ 
 denotes the minimal number of generators of $I$, then $\overline{I^{\mu}} \subset I^* \subset
\bar I$, where for any ideal $J$,  $\overline{J}$ denotes  the integral closure of $J$.
\item[\em (5)] If $R \rightarrow S$ is any ring map, $I^*S \subset (IS)^* $ ("Persistence"). 
\end{enumerate}
}

For the remainder of this section, we will discuss these five main properties, their
proofs and main consequences.
 Some of the five require some mild hypotheses; precise statements will be given. 
All of them are true in any equicharacteristic ring (although Property 2 is not
interesting in characteristic zero). All of them are quite elementary to  prove, at
least in the main settings, with the exception of Property 5 which requires a new idea. 
We will stick to the prime characteristic case, and
simply remark  that "by reduction to characteristic $p$", one can prove the characteristic
zero case  without essential difficulty.

\medskip
Note that one
important property is  omitted from the list. Any decent closure operation
ought to commute with localization, but amazingly, we still do not know that tight
closure does.
\par
\bigskip
{\bf Open Problem} {\it If $U$ is a multiplicative system in a
ring $R$, is
$$
I^*R[U^{-1}] = (IR[U^{-1}])^* ?
$$  }

It is easy to see that one direction  holds,  namely, 
$I^*R[U^{-1}] \subset (IR[U^{-1}])^*$.
 Indeed, if $z \in I^*$, then we have the equations
$cz^{p^e} = a_{1e}y_1^{p^e} + \dots + a_{re}y_r^{p^e}$ in $R$. Expanding to $R[U^{-1}]$,
the same equations show that $\frac {z}{1}$ is in the tight closure of $IR[U^{-1}]$.
This  is a very special case of Property 5 above.
On the other hand, the other direction is not known in any non-trivial case (see, however,
\cite{AHH}, \cite{S6}). 

The localization problem is
 probably the biggest open problem in tight closure theory.
It is remarkable that the theory is so powerful while such a basic question remains
unsolved. The power is derived from the five main properties above, which we now discuss.

\bigskip
\noindent
{\bf{Property One:  All ideals are tightly closed in a regular ring.}}
\medskip

It is easy  to see why all ideals are tightly closed in a regular ring. For example,
consider the special case where $(R, m)$ is local  domain and the Frobenius map is finite.
This is not a very restrictive assumption from our point of view, because we are usually
interested in the local case anyway; also  the Frobenius map is finite in a large class of
interesting rings--- for example,
 for any algebra essentially of finite type over a perfect field or for
any complete local ring with a perfect residue field.

We have a descending chain of subrings of $R$
$$
R \supset R^p \supset R^{p^2} \supset R^{p^3} \supset \dots
$$
Because $R$ is regular, the ring $R$ is a free module considered over each one of the
subrings
$R^{p^e}$. Indeed, the Frobenius map is flat for any regular ring,
 but because we have assumed that $R$ is local and  the Frobenius map is
 finite,
we actually get that $R$ is free over $R^{p^e}$.
This means that, for any non-zero $c$,  we  can find an $R^{p^e}$-linear splitting 
$$
R \stackrel{\phi}{\rightarrow} R^{p^e}
$$
$$
c \mapsto  1
$$
so long as $e$ is large enough that $c$ is not in the expansion of the maximal ideal of
$R^{p^e}$ to $R$
\,(that is,  $c \notin m^{[p^e]}$, where $m^{[p^e]}$ denotes the ideal of $R$ generated by the
$p^{e}-th$ powers of the generators of $m$).

Now if we have an ideal $I = (y_1, \dots, y_r)$ of $R$ and an element $z \in I^*$,
then we can find equations
$$
cz^{p^e} = a_1y_1^{p^e} + \dots + a_ry_r^{p^e}
$$
for all large $e$.
Applying the $R^{p^e}$-linear map $\phi$ above, we see that
$$ z^{p^e} = \phi(a_1)y^{p^e} + \dots + \phi(a_r) y_r^{p^e}$$
where each coefficient $\phi(a_i)$ is some element of $R^{p^e}$.
By taking the $p^{e}-th$ root of this equation,
we see that $z$ is an $R$-linear combination of $y_1, \dots, y_r$.
Thus $z \in I$, and $I^* = I$ for all ideals of $R$.

This completes the proof that all ideals are tightly closed in a regular ring,
at least in the special case we considered. The general case (of prime characteristic)
is not much harder. The point is that flatness of Frobenius in a regular ring. See
\cite{HH1}. 

\bigskip
\noindent
{\bf{Property Two: Elements mapped to $I$ after integral extension are in $I^*$.}}
\medskip

We now prove Property 2: if $R \hookrightarrow S$ is an integral extension of domains of
prime characteristic, and $I$ is an ideal of $R$, then $IS\cap R \subset I^*$.

The following lemma will be useful also in the proof of Property 3.

\medskip
\noindent
{\bf Key Lemma} {\it If $R \hookrightarrow S$ is a module finite extension of domains,
and $d$ is any fixed non--zero element of $S$, then there is an $R$-linear map,
$S \stackrel{\phi}{\rightarrow} R$ sending $d$ to a non-zero element of $R$.}
\medskip

The point in the proof of the Lemma is that after tensoring with the fraction field, $K$,
of $R$, we have an inclusion $ K \hookrightarrow K \otimes_R S$, where
the latter  is simply a finite
dimensional vector  space over $K$. So of course there is a $K$-linear splitting 
$K \otimes S \stackrel{\psi}{\rightarrow} K$ sending
$1 \otimes d$ to $1$. Thinking of $S$ as a subset of $K \otimes S$, 
we look at where  $\psi $ sends each of a set of
$R$-module generators  $\{s_1, \dots,  s_d\}$ for
$S$, say $\psi(s_i) = \frac{r_i}{t_i} \in K$. Now we can define $\phi$ to be
the map $t\psi$ where $t $ is the product of the $t_i$. This map is $R$-linear,
and  sends each $s_i$ to an element of $R$.
The  lemma is proved.

\medskip

To prove Property 2,
let $z \in R$ be any element in $IS\cap R$. This means we can write
$$
z = a_1y_1 + \dots + a_ry_r 
$$
where $a_i \in S$ and the $y_i$'s generate $I$. Because this expression involves only
finitely many
elements from $S$ there is  no loss of generality in assuming $S$ is module
finite over $R$.
Now, raising this equation to the
$p^e-th$ power, we have
$$
z^{p^e} = a_1^{p^e}y_1^{p^e} + \dots + a_r^{p^e}y_r^{p^e}. 
$$
Using the lemma, we find an $R$-linear map $S \stackrel{\phi}{\rightarrow} R$
sending  1 to some non-zero element $c \in R$.   Applying
$\phi$  to this equation, we have
$$
cz^{p^e} = \phi(a_1^{p^e})y_1^{p^e} + \dots + \phi(a_r^{p^e})y_r^{p^e}. 
$$
This is an equation now in $R$, showing that $z \in I^*$. Property 2 is proved.

\medskip
Essentially the same argument shows the stronger property: if $R \hookrightarrow S$
is an integral extension of prime characteristic domains and $I$ is an ideal of $R$, then 
$(IS)^* \cap R \subset
I^*$.
\bigskip

Property 2, unlike the other four properties, is interesting only  in prime characteristic.
For example, if $R$ is any normal domain containing $\QQ$, then $R$ splits off of every 
finite integral extension $S$ (using the trace map).
In this case, $IS\cap R = I$ for every ideal of $R$ and every  integral extension $S$. 

\medskip 
Property 2 can be phrased in terms of the absolute integral closure.
For any domain $R$, the absolute integral closure $R^+$ of $R$ is 
the integral closure of $R$ in an algebraic closure of its fraction field.
In other words, $R^+$ is the direct limit of all finite integral extensions of $R$.
 Property 2  can be stated:
 $IR^+ \cap R \subset I^*$ for all ideals $I$ of $R$. This
leads to the following interesting problem.
\par
\bigskip
{\bf Open Problem}
{\it Let $R$ be a domain of prime characteristic.
Is $IR^+\cap R = I^*$ for all ideals $R$? }

\medskip
In addition to providing a very nice characterization of tight closure, an affirmative answer
to this question would immediately solve the localization problem. Indeed, it is easy to check
that the closure operation defined by expansion to the absolute integral closure and
contraction back to $R$ commutes with localization.   

There is  no non-trivial class of rings in which this open problem has been solved.
However, we do have the following result.

\begin{theorem} {\em \cite{S1}}
Let $R$ be a locally excellent{\footnote{Virtually all rings the commutative algebraist on
the street is likely to run across are locally excellent, but see \cite{Mats} for a
definition.}} domain of prime characteristic. Then
$I^* = IR^+\cap R$ for all parameter ideals $I$ of $R$.
\end{theorem}

A  ''parameter ideal" is  any
ideal $I$ generated by $n$-elements where $n$ is the height of $I$; if $R$ is local, an ideal
is a parameter ideal if and only if it is generated by part of a system of parameters.

\smallskip 
As we see from the theorem,
 tight closure commutes with localization for parameter ideals. However, 
this does not follow from the theorem because this fact  is used in its proof. See
instead
\cite{AHH}.

The proof of this theorem is somewhat involved, so we do not  sketch it here;
see \cite {S1}. The result has been generalized to a larger class of ideals, including
ideals generated monomials in the parameters,
by Aberbach \cite{A}.

\bigskip
\noindent
{ \bf{Property Three: Colon Capturing.}}
\medskip

Property 3, the colon capturing property of tight closure, is
particularly instrumental in applications of tight closure to problems about
Cohen-Macaulayness. Of course, if $R$ is a Cohen-Macaulay local ring with system 
of parameters $x_1, \dots, x_d$, then by definition,
$$
(x_1, \dots, x_i):x_{i+1} \subset (x_1, \dots, x_i)
$$
for each $i = 1, 2, \dots, d-1$.
Colon capturing says that, even for rings that are  not Cohen-Macaulay,
the colon ideal $
(x_1, \dots, x_i):x_{i+1}$ is at least contained in 
 $ (x_1, \dots, x_i)^*$.   Loosely speaking, tight closure
captures the failure of a ring to be Cohen-Macaulay.

We now prove the colon capturing property of tight closure: if $R$ is a local
domain (satisfying some mild hypothesis to be made soon precise)
and $x_1, \dots, x_d$ is a system of parameters for $R$, then 
$$
(x_1, \dots, x_i):x_{i+1} \subset (x_1, \dots, x_i)^*
$$
for each $i = 1, \dots, d-1$.

\medskip
Let us first  assume that $R$ is complete. In this case, we can express $R$ as a module
finite extension of the power series subring $k[[x_1, \dots, x_d]]$, where $k$ is a field
isomorphic to the residue field of $R$. 

Suppose that $z \in (x_1, \dots, x_i):_Rx_{i+1}$.
Consider the ring $A$ contained in $R$ obtained by adjoining the element $z$ to the
power series ring $k[[x_1, \dots, x_d]]$. Observe that the ring $A$ is Cohen-Macaulay;
in fact, $A$ is a hypersurface ring because its dimension is $d$ and its embedding dimension
is $d+1$ (or $d$, if $z$ happens to be in power series ring already).

Now,  because 
$z \in (x_1, \dots, x_i):_Rx_{i+1}$, we can write 
$$z x_{i+1} = a_1x_1  + \dots + a_i x_i$$
for some elements $a_i$ in $R$. Raising this equation to the $p^{e}-th$ power,
we have 
$$z^{p^e} x_{i+1}^{p^e} = a_1^{p^e}x_1^{p^e}  + \dots + a_i^{p^e} x_i^{p^e}$$
Because the inclusion $A \hookrightarrow R$ is a module finite extension,
we can use the Key Lemma to find an $A$--linear map
$R \stackrel{\phi}{\rightarrow} A$ sending $1$ to some
non-zero element $c \in A$. This yields equations
$$cz^{p^e} x_{i+1}^{p^e} = 
\phi(a_1^{p^e})x_1^{p^e}  + \dots + \phi(a_i^{p^e}) x_i^{p^e}$$
where the $\phi(a_j^{p^e})$ are just some elements of $A$.

In other words,
$$
cz^{p^e} \in (x_1^{p^e}, \dots, x_i^{p^e}):_Ax^{p^e}_{i+1}
$$
in the ring $A$. But  $A$ is Cohen-Macaulay, and $x_1^{p^e}, \dots, x_{d}^{p^e}$
is a
system of parameters for $A$, so we see
$$
cz^{p^e} \in (x_1^{p^e}, \dots, x_i^{p^e})
$$
for all $e$.  This shows that $z \in (x_1, \dots, x_i)^* $ in $R$ (also in $A$, 
but  it is $R$ we care about).  Thus 
$ (x_1, \dots, x_i):_Rx_{i+1} \subset (x_1, \dots, x_i)^*,$ and the proof of the colon
capturing property is  complete--- at least for complete local domains.
\bigskip

Inspecting the proof, we see that we have not used  the completeness of $R$ in a crucial
way: what we need is that $R$ the domain is a finite extension of a regular
subring. So this proof also works for algebras essentially of finite type over a field
(the required regular subring is supplied by Noether normalization)
and in many other settings. In fact, colon capturing holds for any ring module finite
and torsion free over a regular ring.
See \cite{HH1} and \cite{Hu2} for different proofs and more general statements.

\medskip
The philosophy of 
colon capturing  holds for other ideals involving  parameters. For example, if $I$ and $J$
are any ideals generated by  monomials in a system of parameters $\{x_0, \dots, x_d\}$,
 one can compute $I:J$
formally as if the $x_i$'s are the indeterminates of a polynomial ring.  Then 
the actual colon $I:J$ is contained in the  tight closure of the 'formal' colon ideal.
Furthermore,  even more is true: we have $I^*:J $ is contained in the
tight closure of the formal colon ideal. Essentially the same proof gives these stronger
 results with very small effort.
For an explicit example, let  $x_0, \dots, x_d$ be a system of parameters in a domain
$R$. Then 
$$
(x_0^t, \dots, x_d^t): (x_0x_1\dots x_d) \subset  (x_0^{t-1}, \dots, x_d^{t-1})^*
$$
and even 
$$
(x_0^t, \dots, x_d^t)^*: (x_0x_1\dots x_d) \subset  (x_0^{t-1}, \dots, x_d^{t-1})^*.
$$
One reason for tight closure's effectiveness is that these sorts of manipulations can often
help
us prove a  general statement  about parameters if we already have an argument for a 
 regular sequence.

\bigskip
\noindent
{\bf{Some Consequences of the First Three Properties.}}

\medskip
It follows immediately from  the colon capturing property that if $R$ is a local
ring in which all ideals
are tightly closed, then $R$ must be Cohen-Macaulay.
Indeed,  if  all parameter ideals
 are
tightly closed,  then colon capturing implies that $R$ is Cohen--Macaulay. 
This leads us to define two important new classes of rings.

\medskip
\noindent
\begin{definition} {\rm A ring $R$ is weakly F-regular if all ideals are tightly closed.
A ring $R$ is F-rational if all parameter  ideals are tightly closed.}
\end{definition}

\medskip
So far we have seen the following implications: Regular $\Longrightarrow$ weakly F-regular 
$\Longrightarrow$ F-rational $\Longrightarrow$ Cohen-Macaulay.
The first implication is Property 1, while the last implication is Property 3.

The reason the adjective "weakly" modifies "F--regular" above goes back to the 
localization problem. 
Unfortunately, we do not know whether the property that all ideals are tightly closed
is preserved under localization. The term "F--regular" is reserved for
rings $R$ in which all ideals are tightly closed  not just in $R$, but also in every
localization of
$R$.   That is, we have the following special case of the localization problem:
\par
\bigskip
\noindent
{\bf Open Problem} {\it If $R$ is weakly F--regular, and $U \subset R$ is any multiplicative
system, is the localization $R[U^{-1}]$ also weakly F--regular?}

\medskip
This problem is  much easier 
than the localization problem itself. Indeed, it has been shown in a number of cases.
For example,  Hochster and Huneke showed the answer is  yes when 
 $R$ is  Gorenstein \cite{HH2}, \cite{HH4}. This was later generalized to the $\Bbb
Q$--Gorenstein case,  and even to the case where there are only isolated non
$\Bbb Q$--Gorenstein points, by MacCrimmon \cite{M}. Using this, it is possible to see that 
weakly F-regular is equivalent to F-regular in dimensions three and less. (These statements
require some mild assumption on $R$, such as excellence). Recently, an affirmative answer was
given also for finitely generated $\NN$-graded algebras over a  field
\cite{LyS}. By contrast, the full
 localization problem  has not been solved in any  of these
cases. 

\medskip
The problem is reminiscent of an analogous problem in commutative algebra that looked quite 
difficult in the mid-century: is the localization of a regular ring still regular?
With Serre's introduction of homological algebra to commutative algebra, the problem
suddenly became quite easy. Perhaps a similar revelation is necessary in tight closure
theory.

\bigskip
Returning to the applications of the first three properties, we 
now prove the following easy, but important, theorem.

\begin{theorem} {\em \cite{HH1}} Let $R \subset S$ be an inclusion of rings that splits
as a map of $R$-modules. If  $S$ is (weakly) F-regular, then $R$ is (weakly) F-regular.
\end{theorem}

The proof is simple.  Suppose that $I$ is an ideal of $R$ and that $z \in I^*$.
This means that for some non-zero $c$, $cz^{p^e} \in I^{[p^e]}$ where $I^{[p^e]}$
denotes the ideal generated by the $p^{e}-th$ powers of the generators of $I$.
Expanding to $S$,  we have
 $cz^{p^e} \in (IS)^{[p^e]}$, so that $z \in (IS)^*$. But all ideals of $S$ are tightly
closed, and so $z \in IS$. Now applying the splitting $S \rightarrow R$ 
(which sends $1$ to $1$ $R$-linearly), we  see that $z \in I$ in $R$ as well.
This completes the proof.

\medskip

The importance of this Theorem lies in the following corollaries.

\begin{corollary} Any ring (containing a field) 
which is a  direct summand of a regular
ring is Cohen-Macaulay.
\end{corollary}

The proof is obvious: a regular ring is F-regular by Property 1, so any direct summand
is also F-regular. By Property 3, this summand is Cohen-Macaulay.

\medskip

\begin{corollary} {\em ({\bf The Hochster-Roberts Theorem})}
The ring of invariants of a linearly reductive  group acting  linearly  on a regular ring
is Cohen-Macaulay.
\end{corollary}

This is essentially a special case of the previous corollary because the
  so-called Reynold's operator
gives us a splitting of $R^G$ from $R$. 

\medskip
We emphasize that both the Theorem and its corollaries make sense and are
true in characteristic zero. Thus even though there are very fewer linearly
reductive groups in prime characteristic, the Hochster-Roberts Theorem for
reductive groups over the complex numbers has been proved  here
by reduction to
characteristic
$p$. To be fair, we have not proved Properties 1 and 3 in characteristic zero (nor
even given a precise definition of tight closure in characteristic zero).
However, if one accepts the existence of a closure operation in characteristic zero
satisfying Properties 1 and 3, then  we have proved that the Hochster-Roberts Theorem follows.

\bigskip
We now mention one of the crown jewels of tight closure theory.
\begin{theorem} {\em \cite{HH3}} Let $R$ be an excellent local domain of
prime characteristic. Then the absolute integral
closure
$R^+$ of $R$  is a Cohen--Macaulay $R$-module.
\end{theorem}

We can see that this must be true as follows. 
Let $x_1, \dots, x_d$  be a system of parameters.
 Suppose $z \in (x_1, \dots, x_i):x_{i+1}$. By  the colon capturing
property, $z \in (x_1, \dots, x_i)^*$. But for parameter ideals, tight closure is the same as
the contraction of the expansion to
$R^+$ (see the  discussion  of Property 2). Thus  $z \in (x_1, \dots,
x_i)R^+\cap R$. This holds for all
$i$, so
 $x_1, \dots,
x_d$ is a regular sequence on $R^+$, and the Theorem is "proved".
Unfortunately, this is not an honest proof  because  the proof that $I^* = IR^+\cap R$
for parameter ideals $I$ uses the Cohen-Macaulayness of $R^+$.

\bigskip
\noindent
{\bf{Property 4: Relationship to integral closure.}}
\medskip

Property 4 is really two statements. First, the tight closure is contained
in the integral closure for any ideal $I$. Second, the integral closure of $I^{\mu}$
(where $\mu$ is the least number of generators of $I$) is contained in the tight closure
$I^*$.

\medskip
The point in proving both statements is the following alternative definition of the integral
closure
$\bar J$ of an ideal
$J$ in a domain $R$: an element $z \in \bar J$ if and only if there exists a non-zero $c$ in
$R$ such that $cz^n \in J^n$  for all (equivalently, for infinitely many) $n \gg 0$.
(This can be easily proved equivalent to  the more standard definition of integral
closure by  recalling another characterization of integral closure:  $\bar J$ consists of 
 all elements $z$ such that $z \in JV$ for all  discrete valuation rings $V$
lying between $R$ and its fraction field. )  Note that in particular, $I^* = \overline{I}$
for all principal ideals $I$. 

Now, with this definition of the integral closure, it is immediately clear that the
tight closure  of any ideal is contained in the integral closure. Indeed, since
the $p^e-th$ power of the generators of $I$ are contained
in the $p^e-th$ power of $I$, we have
$$cz^{p^e} \in I^{[p^e]} \subset I^{p^e}$$ for all $e$.
So any $z $ in $I^*$ is in $\bar I$.

\smallskip
For the second statement, suppose that $ z \in \overline{I^{\mu}}.$
This means that  there exists a non--zero $c$ such that for all $n$,
$cz^n \in I^{\mu n}$.
If $y_1, \dots, y_{\mu}$  generate $I$, then $I^{\mu n}$ is generated by 
monomials of degree $\mu n$ in the $y_i$. But if $y_1^{a_1}y_2^{a_2}\dots y_{\mu}^{a_{\mu}}$
is such a monomial, at least one $a_i$ must be greater than or equal to $n$.
So
$$cz^n \in I^{\mu n} \subset (y_1^n, \dots, y_{\mu}^n)$$
for all $n$.
In particular, this holds for $n = p^e$, for all $e$, and we conclude that $z \in I^*$.
The proof that $\overline{I^{\mu}} \subset I^*$ is complete.

\medskip

The statement that $\overline{I^{\mu}} \subset I^*$  is sometimes called the
Brian\c con-Skoda Theorem.   The original Brian\c con-Skoda Theorem 
stated that for any ideal $I$ in a ring of  convergent complex power series,
the integral closure of the $\mu$-th power of $I$ is contained in $I$,
where $\mu$ is the minimal number of generators of $I$ \cite{BS}.
This statement was later
generalized by Lipman and Sathaye to more general regular local 
rings and then  by Lipman and Tessier to certain ideals in the the so-called 'pseudo-rational' local rings
(for a ring essentially of finite type over a field of characteristic zero, pseudo-rational
is equivalent to rational singularities)
\cite{LS}, \cite{LT}. Tight closure gives an extremely simple proof
of the Brian\c con-Skoda theorem for any regular ring containing a
field: $\overline {I^{\mu}} \subset I^*  \subset
I$, where the first inclusion follows from Property 4 and the second by Property 1.
But better still, tight closure explains what happens in a non-regular ring as well.

The original motivating problem for the Brian\c con-Skoda theorem is said to be due to
J. Mather: if $f$ is a germ of an analytic function vanishing at the origin in $\Bbb C^n$,
find a uniform
$k$ (depending only on $n$) such that $f^k$ is in the ideal generated by the partial
derivatives of $f$.  The Brian\c con-Skoda theorem 
tell us that we can take $k = n$. Indeed, it is easy to check that $f \in 
\overline{ J_f} = \overline{(\frac{\partial f}{\partial x_1}, \dots, \frac{\partial
f}{\partial x_n})}$. So $f^n \in  \overline {J_f}^n \subset \overline {J_f^n} \subset J_f.$ 
It
is remarkable how easy the tight closure proof is for this problem that once seemed very
difficult.

\medskip
Before moving on to Property 5, we consider one more comparison of tight and integral
closure. 
Let $I$ be an $m$-primary ideal in a local domain of dimension $d$. Recall the
Hilbert-Samuel function  defined by $$
HS(n) = \mbox{length} 
R/I^n.
$$
This function 
is eventually a polynomial in $n$, and its normalized leading coefficient
$$
\lim_{n \rightarrow \infty} \frac{d!}{n^d}HS(n)
$$
is called the Hilbert-Samuel multiplicity of $I$.
Analogously, when $R$ is of characteristic $p$, we can define the Hilbert--Kunz function
$$
HK(e) = \mbox{length}
R/I^{[p^e]}.
$$
This function  has polynomial growth in $p^e$, and its leading coefficient
$$
\lim_{n \rightarrow \infty} \frac{1}{(p^e)^d}HK(e)
$$
is called the Hilbert-Kunz multiplicity of $I$.

As is well known, the integral closure of $I$ is the largest ideal containing $I$
having the same Hilbert-Samuel multiplicity (assuming the completion of $R$ is
equidimensional). What is also fairly straightforward to prove is that the tight closure of
$I$ is the largest ideal containing $I$ having the same Hilbert-Kunz multiplicity 
(assuming the completion of $R$ is reduced and equidimensional) \cite{HH1}.
In this sense, tight closure is a natural analog of integral closure.

\medskip
Hilbert--Kunz functions are interesting and mysterious, with important
connections to tight closure theory and  surprising interactions with number theory. 
Much has been proved about them by Paul Monsky, among others; see, for example, \cite{Mo}.
and to the bibliography of \cite{Hu2} for more references on this topic.

\bigskip
\noindent
{\bf{ Property Five: Persistence of Tight Closure.}}

\medskip
The persistence property states: whenever $R \rightarrow S$ is a map of rings containing a field, 
$I^*S \subset (IS)^*$. In other words, any element in the tight closure of an ideal $I$ of
$R$ will ''persist" in being in the tight closure of $I$ after expansion to any $R$--algebra.

Before discussing the precise hypothesis necessary,
let us consider what is involved in proving such a statement.
Suppose $z \in I^*$ where $I$ is an ideal in domain $R$.  Thus there exists a non-zero
$c$ such
that
$$
cz^{[p^e]} \in I^{[p^e]}$$
for all large $e$. 
Expanding to $S$, of course, the same relationship holds in $S$ (using the same letters to
denote
the images of $c$, $z$, and $I$ in $S$).
This would seem to say that the image of $z$ is in $(IS)^*$, 
which is what we need to show. The problem is that $c$ may be in the kernel of the
map $R \rightarrow S$. Thus we need to find a $c$ that  ''witnesses" $z \in I^*$
but is not in this kernel. 

 Unlike the first four properties,
Property 5 does not follow immediately from the definition. The new idea
we need is the idea of a {\it test element.\/}

\medskip
\noindent
\begin{definition} {\rm An element $c$ in a prime characteristic ring  $R$
is said to be in the {\it test ideal\/} of $R$  if, for all ideals $I$ and all elements
$z \in I^*$, we have $cz^{p^e} \in I^{[p^e]}$ for all $e$.  An element $c$ is a {\it test
element\/} if it is in the test ideal but not in any minimal prime of $R$.}
\end{definition}

\medskip
Note that the definition of the  test ideal requires that $cz^{p^e} \in I^{[p^e]}$ for all
$e$,
not just for all sufficiently large $e$. We could also define the {\it asymptotic test
ideal\/} as above but require only that 
$cz^{p^e} \in I^{[p^e]}$ for $e \gg 0$. An interesting fact is that the aymptotic test ideal
is a $D$-module---that is, it is a submodule of the module $R$ under the action of the ring
of all $\ZZ$-linear differential operators on $R$. See \cite{S2}. 
\medskip

It is not at all obvious that there exists a non-zero test ideal 
for a ring $R$. Fortunately, however, it is not very difficult to 
prove the following.

\begin{theorem} {\em \cite{HH2}} Let $R$ be a ring of prime characteristic,
and assume that the Frobenius map of $R$ is
finite. If $c$ is an element of $R$ such that the localization $R_c$ is regular,
then $c$ has a power which is a test element. That is, the test ideal contains
an ideal defining the non--regular locus of  $\spec R$.
\end{theorem}

In a later paper, Hochster and Huneke prove this without the assumption that the Frobenius map
is finite, imposing the weaker  and more technical hypothesis of being finitely generated over
an excellent local ring.  Although the theorem stated above for rings in which 
Frobenius is finite is quite easy to prove, the  proof in the more general setting 
is  difficult and technical; see \cite{HH4}.
\medskip

Note that in any ring $R$, the element 1 is a test element if and only if $R$ is weakly
F-regular. We expect that much more is true:

\medskip
\noindent
{\bf Conjecture.}
{\it The test ideal defines precisely the non--F--regular locus  in $\spec R$.}

\medskip

The conjecture is proved  in some cases, such as for (excellent local) Gorenstein rings
\cite{HH4} and for rings 
$\NN$-graded over a field \cite{LyS}.

\bigskip
Having introduced the idea of a test element, we resume our discussion of
 persistence.
First of all, we should say that Property 5 is not known to hold in the generality
we've stated; some mild hypothesis on $R$  is needed. The problem
is in finding test elements for $R$.

Let us now sketch the proof of persistence.
Let $R \stackrel{\phi}{\rightarrow} S$ be a map of domains.{\footnote{
 Like most proofs in tight closure theory, the proof reduces immediately to the
case where both $R$ and $S$ are domains.}} 
As we remarked above, persistence is trivially true when  $\phi $ is injective,
so  factoring $\phi$ as a surjection followed by an injection,
we might as well assume $\phi$
is surjective.  Now factor $\phi$ as a sequence of surjections
$$
R \rightarrow R/P_1 \rightarrow R/P_2  \rightarrow \dots \rightarrow
 R/(ker \phi) = S,
$$
where $P_1 \subset P_2 \subset \dots \subset (\ker \phi)$
is a saturated chain of prime ideals contained in the kernel of $\phi$. By considering
each map separately, 
we see that we might as well assume that the kernel of the map
$R \stackrel{\phi}{\rightarrow} S$ has height one.

Now if $R$ is normal, then the non-regular locus of $R$ is defined by an ideal
$J$ of height two or more. But as we mentioned above, this means that the test ideal has
height two or more,
{\footnote{This requires some hypothesis on $R$,  such as finite generation over an excellent
local ring,  so that
$R$ satifies the conclusion of Hochster and Huneke's
theorem about test elements above.  In practice, all rings we run across will satisfy this
hypothesis.}} so we can find a
$c$ which is a test element but not in the kernel of $\phi$. The proof is complete in the case
$R$ is normal.

Finally, 
it is not difficult to reduce the problem to the case where $R$ is normal, using Property 2.
What happens is the normalization  $\tilde R$ of $R$ maps to an integral extension $\tilde S$
of $S$,  namely the domain $\tilde S$ obtained by killing a prime of $\tilde R$ lying over
the kernel of $\phi$.  The map $\tilde R \stackrel{\tilde \phi}{\rightarrow}
\tilde S$ restricts
 to the
map $R \stackrel{\phi}{\rightarrow} S$.
  Now if $z \in I^*$ in $R$, then $z \in (I\tilde R)^*$, and
so
$ \tilde \phi(z) 
\in (I\tilde S)^*$ because we know persistence holds when the source ring is normal.
By Property 2 (or really, by the same proof used to prove Property 2), we see that $\phi(z)
\in (I\tilde S)^* \cap S \subset  (IS)^*$. This completes the proof of persistence.

\bigskip

We have completed the proofs and discussion of the five main properties of
tight closure.
It is natural to ask whether  the five main properties characterize tight closure.
They do not, or at least,  not obviously. For example,  in characteristic $p$,
the 'plus
closure'
$IR^+\cap R$ satisfies Properties 1,2, 3, and 5, and in all cases where 
it can be checked, it satisfies Property 4 as well. 
On the other hand, since we expect $I^* = IR^+\cap R$, this is perhaps not very convincing.

A more interesting question is whether we can define a closure operation 
 for rings that do not contain a field (that is, in  'mixed
characteristic') which satisfies Properties 1 through 5. If so, many theorems that can
now be proved only for rings containing a field, such as the homological conjectures
that grew out of Serre's work on multiplicities, would
suddenly admit "tight closure" proofs. The only serious attempts at defining such a closure
operation in mixed characteristic are due to Mel Hochster, but so far none has proved
successful; see, for example, \cite{Ho2}.

\medskip

I hope it is clear from section one that the main ideas in tight closure theory
are remarkably simple and elegant, and also that they have far-reaching consequences.
In section two, we will look more closely at  applications
of tight closure to algebraic geometry.

\bigskip

\sezione{THREE APPLICATIONS OF TIGHT CLOSURE}

\bigskip

At  the  beginning of the part one, we mentioned that tight
closure is applicable to a wide range of problems in commutative algebra and related fields.
Now we will discuss in greater detail how tight closure has
increased our insight in three areas of algebraic geometry:  
adjoint linear systems (Fujita's Freeness Conjecture), vanishing theorems
for cohomology (Kodaira Vanishing), and  singularities.
We will mainly discuss the first of these, giving a tight closure  proof of Fujita's freeness
conjecture for globally generated line bundles, but we point out connections with the other
two topics as they arise.

In all three areas, characteristic $p$ methods are used to prove characteristic zero theorems.
The unifying theme  for the tight closure approach to these three problems is
the action of the Frobenius operator on local cohomology. 

\bigskip
\noindent
{\bf{Reduction to Characteristic $p$.}}

\medskip
Reduction 
to characteristic $p$ is easiest to understand by example.
Say we want to study the affine scheme associated to the ring 
$$
\frac{\Bbb Q[x, y, z]}{(x^3 + y^3 + z^3)}.
$$
We instead consider the "fibration"
$$
\spec \frac{\Bbb Z[x, y, z]}{(x^3 + y^3 + z^3)} \rightarrow \spec  \Bbb Z.
$$
The fiber over a closed point $(p) \in \spec  \Bbb Z$ is the characteristic $p$
scheme
$$\spec \frac{\Bbb Z/(p)[x, y, z]}{(x^3 + y^3 + z^3)},$$
whereas the fiber over the generic point $(0) \in \spec \Bbb Z$
is the original scheme
$$
 \spec \frac{\Bbb Q[x, y, z]}{(x^3 + y^3 + z^3)}.
$$
For the sorts of questions we are interested in here (which are ultimately
cohomological) the following philosophy holds: 
what is true for the generic fiber is true for a Zariski dense set of closed fibers,
and conversely, what is true for a Zariski  dense set of closed fibers is true for the
generic fiber.
 So in order to study  the ring
${\frac{\Bbb Q[x, y, z]}{(x^3 + y^3 + z^3)}}$, we consider the
ring ${\frac{\Bbb Z/(p)[x, y, z]}{(x^3 + y^3 + z^3)}}$ for a "generic" $p$.

\smallskip
The same approach works even if we take $\Bbb  C$ as the ground field.
 Indeed, 
$$
\spec \frac{\Bbb C[x, y, z]}{(x^3 + y^3 + z^3)},
$$
is  obtained from 
$
\spec \frac{\Bbb Q[x, y, z]}{(x^3 + y^3 + z^3)} 
$
by the flat base change $\Bbb Q \rightarrow \Bbb C$.
Again, from the point of view of the types of questions we will consider,
we might as well study $
\frac{\Bbb Q[x, y, z]}{(x^3 + y^3 + z^3)}
$, and hence 
$\frac{\Bbb Z/(p)[x, y, z]}{(x^3 + y^3 + z^3)}$ for a "generic" $p$.

\smallskip
The philosophy holds for any scheme of finite type over a field of
characteristic zero. For example, if we are interested in
the ring
$$
R = \frac{\Bbb C[x, y, z]}{(\pi x^3 + \sqrt{17}y^3 + z^3)},
$$
we  set $A = \Bbb Z[\pi, \sqrt{17}]$ and consider the fibration
$$
\spec \frac{A[x, y, z]}{(\pi x^3 + \sqrt{17} y^3 + z^3)} \rightarrow {\text{ Spec }} A.
$$
Again, we might as well study the prime characteristic 
 ring ${\frac{A/\mu[x, y, z]}{(\pi x^3 +
\sqrt{17} y^3 + z^3)}}$, where $\mu$ is a generic maximal ideal in $A$. Each 
$A/\mu$ is a finite field, so 
these closed fibers are all ''characteristic $p$ models", for varying $p$, of the original
ring
$R$.
\medskip

In general, if 
$$R = \frac{k[x_1, \dots, x_n]}{(F_1, \dots, F_r)},$$
where $k$ is a field of characteristic zero,
we let $A = \Bbb Z[$coefficients of the
$F_i$'s$] \subset k$ and set 
$$ R_A = \frac{A[x_1, \dots, x_n]}{(F_1, \dots, F_r)}.$$
Then the map 
$$
\spec R_A \rightarrow \spec A$$
(or the map $A \rightarrow R_A$)
will be called a  {\it family of models\/} for $\spec R$ (or $R$).
The generic fiber is the original scheme $\spec R$ (after extending the field if necessary)
and a generic (or typical) closed fiber is a {\it characteristic $p$ model\/} of
$\spec R$.
We will prove theorems about $R$ by establishing the same statement
 for a generic characteristic
$p$ model of $R$, that is, "for all large $p$."

\medskip
The idea of a family of models can be used to define concepts in characteristic zero
which seemingly only make sense in prime characteristic. For example, we can define
F-regularity and F-rationality for finitely generated algebras over a field in this way.

\medskip
\noindent
\begin{definition} {\rm Let $R$ be a finitely generated algebra over a field of
characteristic zero. Then $R$ is said to have F-regular type if $R$ admits
a family of models $A \rightarrow R_A $ in which a Zariski dense set of closed fibers are F-regular.
(This does not depend on the choice of the the family of models.) }
\end{definition}

\medskip
Similarly, we can define weakly F-regular type, F-rational type, or F-split type
for any finitely generated algebra over a field of characteristic zero. (In characteristic
$p$, F-split means that the Frobenius map splits, that is, $R^p \subset R$ splits as a
$R^p$-module map.)

There is a subtlety in the meaning of F--regularity for algebras of characteristic zero.
As we've said in part one, the operation of tight closure can be defined for any 
ring containing a field, so it makes sense to define a finitely generated algebra
over a field of characteristic zero  to be weakly F-regular if all ideals are tightly closed.
This is a priori different from the condition of  weakly F-regular type.
We expect that these notions are equivalent, but this remains unsolved. 
See \cite {Ho3}.

\bigskip
The notions of F-rational type  and F--regular type   turn out to be 
intimately connected with the singularities
that come up in the minimal model program. The first theorem in this direction 
explains the name "F-rational".

\begin{theorem} {\em \cite{S3}, \cite{Ha1}}
A finitely generated algebra over a field of characteristic zero has F-rational type
if and only if it has rational singularities.
\end{theorem}

The concept of rational singularities is very important in birational algebraic geometry.
Recall that by
definition, a  ring
$R$ has rational singularities if and only if  it is normal and it admits a desingularization
$X$ for which
$H^i(X, \mathcal O_X) = 0$ for all $i >0$.)
\smallskip

We will not dwell on this theorem here, rather refering the the papers \cite{S3} and 
\cite{Ha1} in the bibliography. 
Later,  we will later mention some ideas in the proof. 
Now we move on the application of tight closure to Fujita's freeness conjecture, where many
related ideas appear.

\bigskip
\noindent
{\bf Application of Tight Closure to Adjoint Linear Series.}
\medskip

Let $X$ be a smooth projective variety of dimension $d$, and let $\mathcal L$ be an 
ample invertible sheaf on $X$.  We are interested in the adjoint line bundles
$\omega_X \otimes \mathcal L^n$, for $n >0$. Because $\mathcal L$ is ample, we know that 
for large $n$, this adjoint bundle is globally generated. Fujita's freeness conjecture
provides an effective version of this statement.

\medskip
\noindent
{\bf Fujita's Freeness Conjecture.}
{\it With $X$ and $\mathcal L$ as above, the sheaf $\omega_X \otimes
\mathcal L^{d+1}$ is globally generated. }

\medskip
The conjecture is known in characteristic zero in dimension up to four \cite{R}, \cite{EL},
\cite{Ka}. See \cite{Ko} for a survey.
 In arbitrary characteristic, the best that is known is
given by the following theorem.

\begin{theorem} {\em \cite{S4}} If $X $ is a smooth projective variety of dimension $d$ and $\mathcal
L$ is a globally generated ample line bundle on $X$, then $\omega_X \otimes \mathcal L^{d+1}$ is
globally generated.
\end{theorem}

See \cite{S7} for a recent improvement of this result.

\bigskip

Our next task is to prove this theorem, that is, to establish Fujita's Freeness Conjecture for 
globally generated line bundles.
This will give a good overview of some of the methods that can be used in applying tight
closure to algebro-geometric questions.

\medskip

If $X$ has characteristic zero, the first step is to reduce to the characteristic $p$ case
using the standard technique we described. So it is enough to prove the theorem 
in the case that $X$ has prime characteristic.

\medskip
A good way to study 
an ample line bundle on a projective variety $X$ is to 
build the section ring 
$$
S = \oplus_{n \geq 0} H^0(X, \mathcal L^n).
$$
This is a finitely generated, $\NN$-graded ring whose associated projective scheme 
recovers $X$. Its dimension is $d+1$. 
Assuming that $X$ is irreducible, every section ring $S$ will be a domain.
If $X$ is smooth, then $S$ has (at worst) an isolated
singularity at the unique homogeneous maximal ideal $m$.  
The invertible sheaf $\mathcal L^n$ corresponds to the graded $S$-module $S(n)$,  the
$S$-module $S$ with degrees shifted by n.

\smallskip
Fujita's Freeness Conjecture is equivalent to the following more commutative algebraic
statement.
\par
\bigskip
\noindent
{\bf Fujita's Freeness Conjecture in terms of local cohomology}
{\it If $(S, m)$ is a section ring with an isolated non-smooth point, 
then $H^{d+1}_m(S)$ has the following property: there exists an integer $N$ 
such that for all $\eta \in H^{d+1}_m(S)$ of degree less than $N$,  $\eta$ has a non-zero
$S$-multiple of degree $-d-1$.}

\medskip
The proof of the equivalence of this statement with Fujita's Conjecture is 
not difficult. This is essentially the dual statement (using Matlis duality for $S$ or 
Serre duality for $X$). Details can be found in  \cite{S4}.

\medskip
To prove Fujita's Conjecture, we will tackle this local cohomological conjecture.
First we describe a convenient way to think about elements in the local cohomology module
$H^{d+1}_m(S)$.

Let  $x_0, x_1, \dots, x_d$ be a system of parameters for $S$ of degree one.
Such a system of parameters exists by our assumption that $\mathcal L$ is globally generated
 (after enlarging the ground field if necessary).
The local cohomology module 
$H^{d+1}_m(S)$ can be computed as the cokernel of the following map

\begin{eqnarray*}
S_{x/x_0} \oplus S_{x/x_1} \oplus \cdots \oplus S_{x/x_n}
\,&\longrightarrow & S_x \\
(\frac{s_0x_0^t}{x^t},\frac{s_1x_1^t}{x^t}, \dots, \frac{s_dx_d^t}{x^t})
 & \mapsto &
\frac{\sum_{i=0}^d (-1)^is_ix_i^t}{x^t}
\end{eqnarray*}
where $x $ denotes the product $x_0x_1\dots x_d$ of the $x_i$'s.   This is the last map in
the Cech complex for computing the cohomology of the sheaf of $\mathcal O_X$-algebras 
$\bigoplus_{n=0}^{\infty} \mathcal O_X(nL)$
 with respect to the affine cover of $X$ given by the $\,d+1\,$ open
sets
$U_i$  where
$x_i$ does not vanish.
More generally, the local cohomology modules 
$H^i_m(S)$ can be computed as the cohomology of this  Cech complex, so that 
$H^i_m(S) = \bigoplus_{n \in \Bbb Z} H^{i-1}(X, \mathcal L^n)$, for all $ i > 1$.

We represent elements of 
 $H^{d+1}_m(S)$  by fractions $[\frac{z}{x^t}]$ , with the square bracket
reminding us of the  equivalence relation on fractions.  
If the degree of $\eta$ is $-n$, we see that $-n = {\text{ deg }} z - t(d+1)$.

It is easy to see that if $z \in (x_0^t, x_1^t, \dots, x_d^t)$,
then $\eta =
[\frac{z}{x^t}]$ must be zero, by thinking about the image of the map above.
Unfortunately, the
converse is false. However, we have the following interesting observation.

\begin{lemma}
If $\eta = [\frac{z}{x^t}] = 0$, then $z \in (x_0^t, \dots, x_d^t)^*$.
\end{lemma}

The Lemma is easily proved:  if $\eta = 
[\frac{z}{x^t}] = 0$, then this means that  for some integer $s$, we have
$$
[\frac{z}{x^t}] = [\frac{x^sz}{x^{t+s}}] = 0
$$
where now $x^sz \in (x_0^{t+s}, \dots, x_d^{t+s})$.
Thus 
$$
z\in 
(x_0^{t+s}, \dots, x_d^{t+s}): x^s$$
so by colon capturing, $z \in (x_0^{t}, \dots, x_d^{t})^*$. The Lemma is proved.

\bigskip
\noindent
{\bf The Frobenius Action on Local Cohomology.}

\medskip
The Frobenius action on local cohomology is the main idea in the proof
of Fujita's Freeness Conjecture for globally generated line bundles, and in
the proof of the equivalence of rational singularities with F-rationality. It is 
also  the central point in the relationship between tight closure and  the Kodaira Vanishing
theorem.
 The idea of using the Frobenius action on local cohomology to study tight closure
first appears in   the work of Richard Fedder and Kei--ichi Watanabe \cite{FW}.

\smallskip
The Frobenius action
$$
H^{d+1}_m(S) \stackrel{F}{\rightarrow} H^{d+1}_m(S) 
$$
is easy to understand. Indeed, Frobenius acts in a natural way  on
each
 module $S_{x_{i_1}\dots
x_{i_r}}$ in the Cech complex  defining the local cohomology modules; it simply
raises
fractions to their $p-th$ powers. This action obviously   commutes with the boundary maps,
so that it induces a natural action on the local cohomology modules.
In particular, the Frobenius action on $H^{d+1}_m(S)$ is given by
$$
\eta = [\frac{z}{x^t}] \mapsto \eta^p  = [\frac{z^p}{(x^t)^p}].
$$

Using this, it makes sense to define tight closure for submodules of $H^{d+1}_m(S)$ by
mimicking the definition for ideals. For example, we can define the tight closure of the zero
submodule in
$H^{d+1}_m(S)$:
$$
0^* = \{ \eta \in H^{d+1}_m(S)  \, |\, {\mbox{ there exists }}
 c \neq
0 {\mbox{ with }} c\eta^{p^e} = 0 {\mbox{ for all }} e \gg 0\}.
$$
The tight closure of zero in $H^{d+1}_m(S)$ is an important gadget. One can show that 
it is the unique maximal proper submodule of $H^{d+1}_m(S)$ stable under the action of
Frobenius \cite{S3}.

\bigskip

Returning to the proof of Fujita's Freeness Conjecture,  we observe the
 following two facts.
\begin{enumerate}
\item[(1)] $\eta = [\frac{z}{x^t}] \in 0^\ast$ if and only if
$z \in (x_0^t, \dots, x_d^t)^\ast$.
\item[(2)] Any test element $c$ kills $0^\ast$.
\end{enumerate}
These two facts are straightforward to prove using nearly the same argument  as in the  
proof of the Lemma above.

\smallskip
Now the proof can be summarized in two main steps.
First we show that if $\eta \in H^{d+1}_m(S)$  
does not have a multiple of degree $-d-1$, then $\eta$ is in $0^*$.
Next we show that $0^*$ vanishes in all sufficiently small degrees.
Obviously, upon completion of these two steps, the proof is complete.

\bigskip
{\it Step One: if $\eta \in H^{d+1}_m(S)$  
does not have a multiple of degree $-d-1$, then $\eta$ is in $0^*$. \/}

\medskip The main point is colon capturing. 
Assume on the contrary, that an element 
$\eta = [\frac{z}{x^t}]$ of degree $-n$ has no non-zero multiple of degree $-d-1$.
This means that every element of $S$ of degree $n - d-1$ must kill $\eta$.
In particular, 
$$
(x_0, \dots, x_d)^{n - d -1} \eta = 0.
$$
By the Lemma, this means that 
$$
(x_0, \dots, x_d)^{n - d -1}z \in 
(x_0^t, \dots, x_d^t)^{*},$$
or in other words, 
$$
z \in (x_0^t, \dots, x_d^t)^{*}: (x_0, \dots, x_d)^{n - d -1}.
$$

Now we  use colon capturing. We manipulate the parameters  $x_0, \dots, x_d$ formally as if
they are the indeterminants of a polynomial ring, in which case the colon 
ideal (ignoring the $*$) would be easily computed to be 
$$
(x_0^t, \dots, x_d^t) +  (x_0, \dots, x_d)^{(d+1)(t-1) - (n - d -1) +1}.
$$
Colon capturing  
says that the actual colon ideal is contained in the tight closure of this "formal" colon
ideal, that is, 
$$
z \in [(x_0^t, \dots, x_d^t) +  (x_0, \dots, x_d)^{(d+1)(t) - (n) +1}]^*.
$$
But note that the degree $z$ is $(d+1)t -n $ (because $\eta = 
 [\frac{z}{x^t}]$ has degree $-n = {\text{ deg }} z -(d+1)t$).
Thus
$$
z \in [(x_0^t, \dots, x_d^t) +  (x_0, \dots, x_d)^{\deg z +1}]^*.
$$
A moment's thought reveals that this forces
$$
z \in (x_0^t, \dots, x_d^t)^*.
$$
Indeed, if
$$
cz^q \in (x_0^t, \dots, x_d^t)^{[q]} +  [(x_0, \dots, x_d)^{\deg z +1}]^{[q]},
$$ we see immediately that because the degree of $c$ is fixed, the degrees of the
generators of $[(x_0, \dots, x_d)^{\deg z +1}]^{[q]}$ are much larger than the degree
of $cz^q$, so  that $cz^q$ must in fact be in the ideal $
(x_0^t, \dots, x_d^t)^{[q]} $ for large $q = p^e$.
But by Fact  (1) above, then we see that 
$$
\eta = 
 [\frac{z}{x^t}] \in 0^*,
$$ 
and the proof of step one is complete.

\bigskip
{\it Step two: $0^*$ vanishes in  sufficiently small degrees. }

\medskip
The point is to consider the test elements  of $S$.
Because $X$ is smooth, the section ring $S$ has an isolated singularity.
This means that the defining ideal of the non-regular locus of $S$ is $m$-primary.
As we  mentioned in Lecture 1, this implies that the test ideal  of $S$ (of all elements
that ''witness" all tight closure relations)  contains an  $m$-primary ideal.
But according to Fact 2 above, the test ideal of $S$ annihilates $0^*$,
so that $0^*$ is killed by an $m$-primary ideal. This says that $0^*$ has finite length,
so of course, it must vanish eventually in all degrees sufficiently small.
This completes the proof of step two, and thus the proof of Fujita's Freeness Conjecture
for globally generated line bundles.

\bigskip
Experts will notice that the argument above does not really require that $X$ be smooth.
We used smoothness only in Step 2,   to conclude that $0^*$ is finite length.
But  $0^*$ is of finite length more generally, and is in fact equivalent to the variety 
$X$  being F-rational (or F-rational type in characteristic zero).
Thus Fujita's Freeness Conjecture holds for any globally generated ample line bundle on a
projective F-rational (type) variety.

\bigskip
We should remark that Fujita's Freeness Conjecture for  globally generated line bundles
can also be proved, in characteristic zero, using the Kodaira vanishing theorem.
As far as I know, however, tight closure provides the  only proof in prime characteristic.
Interestingly, 
the Frobenius action on local cohomology  seems to act as a substitute for Kodaira
Vanishing. 
There is a good reason for this:
it turns out that 
Kodaira vanishing theorem is equivalent to a
statement about the action on Frobenius on local cohomology modules.

\bigskip
\noindent
{\bf Tight Closure and  Kodaira Vanishing.}
\medskip

Recall the classical Kodaira Vanishing Theorem:

\medskip
{\bf Kodaira Vanishing}
{\it If $X$ is a smooth projective variety of characteristic zero,
and $\mathcal L$ is any ample invertible sheaf on $X$, then $H^i(X, \mathcal L^{-1}) = 0$
for all $i < \dim X$. }

\medskip
The Kodaira Vanishing Theorem 
is false in characteristic $p$, although it can be proved by reduction to characteristic $p$
 \cite{DI}. See also \cite{EV}.

\smallskip
Let $S = \oplus_{n \geq 0} H^0(X, \mathcal L^n)$ be the section ring of the pair $(X, \mathcal L)$.
Unwinding definitions using the point of
view that local cohomology can be computed from the Cech complex of the $\mathcal O_X$--algebra
$\oplus \mathcal L^n$, Kodaira Vanishing 
is  seen to be equivalent to 
$$
H^i_m(S) {\mbox{ vanishes in negative degree for all }} i {\mbox{ with }}
1 < i < \dim S.
$$

Because $S$  has (at worst) an isolated non-Cohen-Macaulay point at $m$,
we know that each $H^i_m(S)$ is supported at $m$, and hence must vanish
 in degrees sufficiently small. So we could state the 
Kodaira Vanishing Theorem as follows:
{\it the Frobenius action on  a dense set of characteristic $p$ models for $S$ 
is injective  in negative degrees  on $H^i_m(S)$, for $ 1 < i < \dim S$.\/}

Although it may  sound a bit silly, this way of stating the vanishing of local 
cohomology in  negative degree has the advantage of  making sense also for the 
top local cohomology module $H^{\dim S}_m(S)$. In fact, the injective action of
Frobenius  on  $H^{d+1}_m(S)$ in negative degrees is a new  and important phenomenon, a natural
generalization of  the Kodaira Vanishing Theorem, which is not at all apparent otherwise. 
 This extension to the top local
cohomology module  was conjectured to be true and called
''Strong
       Kodaira Vanishing" in \cite{HS}. The conjecture was proved in a beautiful
paper of Nobuo Hara \cite{Ha1}, and in fact, turns out to be the main point in his proof
that a rationally  singular variety (of characteristic zero)  must be of  F--rational type.
(See also \cite{MS}.)

\medskip
The injective action of Frobenius on the negative degree part of local cohomolgy can be 
re-interpreted in terms of tight closure of parameter ideals. Using ideas similar
to the ideas we used in the proof of Fujita's Conjecture to translate statements about 
the Frobenius action on $\eta = [\frac{z}{x^t}]$  into statements about the tight closure
of $(x_0^t, x_1^t, \dots, x_d^t)$, we get a tight closure version of Kodaira Vanishing.

\medskip
\noindent
{\bf Kodaira Vanishing in terms of Tight Closure.} \cite{HS}
{\it Let $S$ be a  section ring of a pair $(X, \mathcal L)$
where $X$ is a smooth variety of characteristic zero and $\mathcal L$ is an ample invertible 
sheaf of $\mathcal O_X$-modules.  Then for  any proper subset $x_0, \dots, x_k$ of a system of
(homogeneous) parameters for $S$, where $\deg x_i \gg 0$,
$$
(x_0, \dots, x_k)^* \subset  \sum_{i = 0}^k (x_0, \dots, \hat x_i, \dots, x_k)^*  + S_{\geq D}
$$
where $D$ is the sum of the degrees of the $x_i$'s.}

This theorem is  equivalent to the Kodaira Vanishing Theorem.
Just as Kodaira Vanishing can fail in prime characteristic, so can this tight closure
statement. However, the statement holds when  $S$ is a generic characteristic $p$ model
for a section ring of characteristic zero, that is, ''for large $p$."

By allowing the possibility that we have a full system of parameters in the above version of
the Kodaira Vanishing Theorem,
we get the strong  Kodaira Vanishing Theorem.  In fact, if $x_0, x_1, \dots, x_d$ is a full
system of parameters for a section $S$ as above, we get a more precise statement.

\medskip
\noindent
{\bf Strong Kodaira Vanishing} \cite{HS} \cite{Ha1}.
{\it Let $S$ be an $\Bbb N$-graded ring over a field of characteristic zero, and let 
$x_0, x_1, \dots, x_d$ be a full system of (homogeneous) parameters for $S$, with 
$\deg x_i \gg 0$.
Then
$$
(x_0, \dots, x_d)^* =  \sum_{i = 0}^d (x_0, \dots, \hat x_i, \dots, x_d)^*  + S_{\geq
D}
$$
where $D$ is the sum of the degrees of the $x_i$'s. }

\medskip
The reason we get equality here is that $S_{\geq D} $ is contained in $(x_0, \dots, x_d)^*,$
as can be verified with the Brian\c con-Skoda theorem (Property 4).

\smallskip

It is possible to say precisely how large the degrees of the $x_i$'s must be in the statements
of Kodaira and strong Kodaira vanishing in terms of tight closure.
In both theorems, each $x_i$ should have degree larger than $a$, where $a$ is the a-invariant
of $S$. By definition (due to Goto and Watanabe), the  $a$-invariant  is the largest 
integer $n$ such that $H^{\dim S}_m(S)$ is non-zero in degree $n$.

The strong form of Kodaira Vanishing  is conjectured in \cite{HS}, where the idea of
the  ''monomial property of a $d^+$ sequence" due to Goto and Yamagishi is used. It is proved
in
\cite{HS} for rings of dimension two, from which  it is shown that
the Kodaira Vanishing Theorem follows for any normal surface of dimension two. 
In full generality, however, the statement was not known until Nobuo Hara proved
the injectivity of the Frobenius action on the negatively graded part of local cohomology
\cite{Ha1}.  Hara has since greatly generalized his work; see \cite{Ha3}.
\bigskip

\noindent
{\bf Tight Closure and  Singularities.}

\medskip
Finally, we summarize some more connections
  between tight
closure and  singularities in algebraic geometry.  

\smallskip
Let $X$ be a normal variety of characteristic zero.
Assume that $X$ is $\Bbb Q$-Gorenstein, that is, that the reflexive sheaf $\omega_X$ 
represents a torsion element  $K_X$ in the (local) class group of $X$. In other words,
the Weil divisor class $K_X$ is assumed to have a multiple which is locally principal.

Consider a desingularization $\tilde X \stackrel{\pi}{\rightarrow} X$ of $X$, where the exceptional
divisor is a simple normal crossings divisor with components $E_1, \dots, E_n$.
Write 
$$
K_{\tilde X} = \pi^*K_X + \sum_{i= 1}^n a_iE_i
$$
for some unique rational numbers $a_i$. To understand this expression,
suppose that $rK_X$ is locally principal, so that it makes sense to pull it back;
then compare to $rK_{\tilde X}$. The difference is some divisor supported on the
exceptional set, hence of the form $\sum_{i= 1}^n m_iE_i$. Dividing by $r$,
we arrive at the above expression, where 'equality' means numerical equivalence of
$\Bbb Q$-divisors. See \cite{KMM}.

\medskip

In general, the $a_i$'s can be any rational number, although if $X$ is  smooth, we can
easily see that each $a_i$ will be a positive integer. This leads us to the following
restricted class of singularities.

\medskip
\noindent
\begin{definition} {\rm The variety $X$ has log-terminal singularities if all $a_i > -1$, and has log-canonical
singularities if all $a_i \geq 1$. (This is independent of the
choice of desingularization.) }
\end{definition}

\medskip

The relationship to tight closure
is is evidenced by the following theorem.

\begin{theorem} Let   $X$ be a normal $\Bbb Q$-Gorenstein variety  of characteristic zero.
$X$ has F-regular type  if and only if $X$ has log-terminal singularities. 
\end{theorem}

This theorem follows immediately from the 
equivalence of rational singularities and F-rational type discussed earlier,
using the "canonical cover trick".
Indeed,  assuming $X$ is local, set  
$$
Y = \spec \{\mathcal O_X \oplus \mathcal O_X(K_X) \oplus \mathcal O_X(2K_X) \oplus \dots
\mathcal O_X((r-1)K_X)\}
$$ where $r$ is such that $\mathcal O_X(rK_X) $ is isomorphic to  $\mathcal O_X$ via a fixed
isomorphism (so that we can define a ring structure on 
$\mathcal O_X \oplus \mathcal O_X(K_X) \oplus \mathcal O_X(2K_X) \oplus \dots
\mathcal O_X((r-1)K_X)$).
The  natural map $Y \rightarrow X$ 
is called the canonical cover of $X$. It is easy to check that when $X$ is Cohen-Macaulay,
the canonical cover $Y$ is Gorenstein, and that the map is \'etale in codimension one.
With these properties, it is   not hard to show the following two facts:
\begin{enumerate}
\item[(1)] (Kawamata)
$Y$ has rational singularities if and only if $X$ has log-terminal singularities.
\item[(2)] (K.-i. Watanabe)
Y has F-rational type if and only if $X$ has F-regular type.
\end{enumerate}
Thus the equivalence of F-regular type with log-terminal singularities follows
from the equivalence of F-rational type with rational singularities.
 \bigskip

There are some
subtleties involved in the argument using the canonical cover.
Watanabe's argument shows  F-rationality for $Y$ 
is equivalent to {\it strong F-regularity\/} for $X$. Strong F-regularity is
a technical condition conjectured to be equivalent to weak F-regularity (when both
are defined), introduced because it, unlike weak F-regularity, is easily shown to  pass to
localizations  \cite{HH2}. However,  in the case of
$\Bbb Q$--Gorenstein rings, weak and strong F-regularity turn out to be equivalent \cite{M}.

\smallskip
The first proof  that F-regular type $\Bbb Q$-Gorenstein singularities are log-terminal 
is due to  Kei-ichi Watanabe and uses a different argument \cite{W}. 
This different argument also produces the following nice result.

\begin{theorem} {\em \cite{W}}
Let $X$ be a variety satisfying the conditions above.  If $X$ is of F--split type,
then $X$ has log--terminal singularities. (Recall, a local ring of characteristic $p$ is
F-split if the inclusion $R^p \hookrightarrow R$ splits as a map of $R^p$ modules.) 
\end{theorem}

A very interesting open problem that has deep connections with number theory is
the following.
\par
\bigskip
\noindent
{\bf Open Problem.}
{\it If $X$ has log-canonical singularities, does $X$ have F--split type?}

\bigskip
\noindent
{\bf Further Reading on Tight Closure.\/}

\smallskip
The original tight closure paper of Hochster and Huneke \cite{HH1} is still an excellent
introduction to the subject.
There are  also a number of expository articles on tight closure.
Craig Huneke's book {\it Tight Closure and its Applications\/} \cite{Hu2}
is an good place for a beginning commutative algebra student to
learn the subject; it contains several applications more or less disjoint from
the ones discussed in detail  here. It also contains an appendix by Mel Hochster \cite{Ho3}
discussing tight closure in characteristic zero. 
Another nice survey is \cite {Ho1}, which contains a list of
open problems; although the article is now seven years old, many of these problems 
remain open. A more recent view is provided by  the expository article 
\cite{B}.
The article \cite{S5} is a survey written for algebraic geometers.
Huneke's ''Tight Closure and Geometry" is another nice read for algebraists
\cite{Hu3}. All these sources, 
but especially \cite{Hu2}, contain long bibliographies to direct the
reader to numerous research articles on tight closure.

\end{document}